# Limit Theory of the Multi-set Allocation Occupancy (MAO) Distribution: Normal and Poisson Approximations via MAO Norm

Xing-gang Mao

Department of Neurosurgery, Xijing Hospital, the Fourth Military Medical University, Xi'an, No. 17 Changle West Road, Xi'an, Shaanxi Province, China;

*Corresponding authors. E-mail addresses: xgmao@fmmu.edu.cn (X.G. Mao);

ORCID: 0000-0001-8505-1904

Short title: Normal and Poisson Approximations of MAO Distribution

Conflict of Interest Disclosures: None reported.

Funding/Support: This work was supported by the National Natural Science Foundation of China (Grant No: 82473137), Xijing Innovation Research Institute - Joint Innovation Fund (LHJJ24YX19).

# Abstract

This paper investigates the asymptotic behavior of the Multi-set Allocation Occupancy (MAO) distribution, which models the count vector $\mathbf{X} = (X_{=0}, \ldots, X_{=T})$ from $T$ independent rounds of sampling without replacement of size $m$ from $N$ individuals. Focusing on $X_{=t}$ (individuals in exactly $t$ subsets) and employing the MAO norm—a combinatorial tool yielding closed-form factorial moments—we derive the exact marginal distribution of a single individual as $\text{Bin}(T, p)$ with $p = m/N$. Using the MAO norm, we prove that for any fixed number of distinct individuals, their joint distribution differs from the product of marginals by $O(1/N)$, strengthening the elementary pairwise covariance estimate and establishing the weak dependence required for limit theorems.

Based on these findings, we delineate two asymptotic regimes:

- **Normal approximation:** When $N \to \infty$ with $p = m/N$ fixed, $X_{=t}$ obeys a central limit theorem and can be approximated by a normal distribution with mean $N\pi_t$ and variance obtained from the MAO norm, where $\pi_t = \binom{T}{t} p^t (1-p)^{T-t}$.

- **Poisson approximation:** When the expected value $\lambda_N := N\pi_t \to 0$, $X_{=t}$ converges in distribution to $\text{Poisson}(\lambda)$. By symmetry, if $\lambda_N \to N$, then $N - X_{=t}$ converges to a Poisson distribution. In the intermediate regime where $\lambda_N/N$ stays away from $0$ and $1$, the normal approximation applies.

All theoretical results are rigorously proved via moment methods and are corroborated by extensive numerical simulations, which further demonstrate that the approximations remain valid even when the subset sizes are not equal. The MAO norm thus emerges as a unifying tool that connects the exact combinatorial structure of the model to its asymptotic theory, while the MAO inequality provides a powerful way to estimate the order of dependence. These findings offer clear guidance for choosing the appropriate approximation in practical applications, and they open the door to using the MAO distribution as a null model for detecting non-random aggregation in fields as diverse as ecology, sociology, and statistical physics.

## 1. Introduction

Consider a finite population of $N$ individuals, labeled $1, 2, \ldots, N$. There are $T$ subsets $S_1, S_2, \ldots, S_T$, each of size $m$ ($0 < m < N$), drawn independently from the population by simple random sampling without replacement. Let $p = m/N$, and assume that $p$ remains fixed as $N \to \infty$. For a fixed integer $t$ ($0 \leq t \leq T$), define the random variable

$$X_{=t} = \#\{i \in \{1, \ldots, N\} : \text{individual } i \text{ belongs to exactly } t \text{ subsets}\}$$

We are interested in the distributional properties of the random vector $\mathbf{X} = (X_{=0}, X_{=1}, \ldots, X_{=T})$. This paper systematically develops the fundamental concepts of the MAO distribution and establishes its asymptotic relationships with the binomial, Poisson, and normal distributions. The present work builds on the foundations laid in previous studies [1-3], which introduced the MAO function, its moment formulas, which actually substantially extended the traditional general hypergeometric distribution and reveals its profound intrinsic rich structure.

## 2. Definition of the MAO Norm

We first introduce the central concept of MAO theory—the MAO norm. Let $[T] = \{1, \ldots, T\}$. For any collection of subsets $A_1, \ldots, A_r \subseteq [T]$, define the function

$$g(A_1, \ldots, A_r) = \prod_{i=1}^{T} (m)_{k_i} (N-m)_{r-k_i}$$

where $(x)_k = x(x-1)\cdots(x-k+1)$ denotes the falling factorial (with $(x)_0 = 1$), and $k_i = \sum_{j=1}^{r} \mathbf{1}_{i \in A_j}$ is the number of times index $i$ appears in these $r$ subsets.

For a given list of nonnegative integer parameters $(p_1, \ldots, p_r)$ (each $p_j \in \{0, 1, \ldots, T\}$), define the transversal sum

$$G_T(p_1, \ldots, p_r) = \sum_{\substack{A_1, \ldots, A_r \subseteq [T] \\ |A_j| = p_j}} g(A_1, \ldots, A_r).$$

The MAO norm is then defined as

$$\| (p_1, \ldots, p_r) \|_T = \frac{G_T(p_1, \ldots, p_r)}{(N)_r^{T-1}}.$$

Here $(N)_r$ is the falling factorial, and $(N)_r^{T-1}$ denotes its $(T-1)$-th power.

## 3. Meaning of the MAO Norm and Moment Formulas

The importance of the MAO norm lies in its direct provision of moments for various random variables in the MAO distribution. In particular, for the random vector $\mathbf{X} = (X_{=0}, \ldots, X_{=T})$, any factorial moment can be expressed as[3]

$$\mathbb{E}\left[\prod_{t=0}^{T}\left((X_{=t})_{n_t}\right)\right] = \| \underbrace{0, \ldots, 0}_{n_0}, \underbrace{1, \ldots, 1}_{n_1}, \ldots, \underbrace{T, \ldots, T}_{n_T} \|_T$$

where $(X)_n = X(X-1)\cdots(X-n+1)$ is the falling factorial, and the right-hand side is the norm obtained by listing each category $t$ repeatedly $n_t$ times. This formula covers all possible product moments and serves as a unified tool for all moment calculations in the MAO distribution. It reveals that all moment information of the MAO distribution is encoded in the norm, providing a unified framework for subsequent asymptotic analysis.

By choosing specific $n_t$, we obtain various statistics:

- **First moment (expectation):** Taking $n_t = 1$ and all others zero gives

$$\mathbb{E}[X_{=t}] = \| t \|_T.$$

- **Second moments and covariance:**

Second factorial moment: $n_t = 2$ yields

$$\mathbb{E}[X_{=t}(X_{=t} - 1)] = \| t, t \|_T.$$

Hence the variance is

$$\mathrm{Var}(X_{=t}) = \| t, t \|_T + \| t \|_T - (\| t \|_T)^2.$$

Product moment of two different variables: $n_t = 1, n_s = 1$ $(t \neq s)$ gives

$$\mathbb{E}[X_{=t} X_{=s}] = \| t, s \|_T,$$

and the covariance is

$$\mathrm{Cov}(X_{=t}, X_{=s}) = \| t, s \|_T - \| t \|_T \| s \|_T.$$

- **Third moments and third central moments:**

Product moment of three variables: $n_t = 1, n_s = 1, n_u = 1$ gives

$$\mathbb{E}[X_{=t} X_{=s} X_{=u}] = \| t, s, u \|_T.$$

Third mixed central moment: Let $\mu_t = \mathbb{E}[X_{=t}]$. Then

$$\mathbb{E}[(X_{=t} - \mu_t)(X_{=s} - \mu_s)(X_{=u} - \mu_u)]$$

$$= \| t, s, u \|_T - \mu_t \| s, u \|_T - \mu_s \| t, u \|_T - \mu_u \| t, s \|_T + 2\mu_t \mu_s \mu_u$$

- **Fourth moments and fourth central moments:**

    Product moment of four variables: $n_t = 1, n_s = 1, n_u = 1, n_v = 1$ gives

    $$\mathbb{E}[X_{=t}X_{=s}X_{=u}X_{=v}] = \| t, s, u, v \|_T.$$

    Fourth mixed central moment: For the case of four identical variables, this yields quantities related to kurtosis. The general form is

    $$\mathbb{E}[(X_{=t} - \mu_t)(X_{=s} - \mu_s)(X_{=u} - \mu_u)(X_{=v} - \mu_v)]$$

    $$\begin{aligned}= \quad & \| t, s, u, v \|_T - \mu_t \| s, u, v \|_T - \mu_s \| t, u, v \|_T \\ & -\mu_u \| t, s, v \|_T - \mu_v \| t, s, u \|_T + \mu_t \mu_s \| u, v \|_T \\ & +\mu_t \mu_u \| s, v \|_T + \mu_t \mu_v \| s, u \|_T + \mu_s \mu_u \| t, v \|_T \\ & +\mu_s \mu_v \| t, u \|_T + \mu_u \mu_v \| t, s \|_T - 3\mu_t \mu_s \mu_u \mu_v.\end{aligned}$$

This expression is obtained by expanding the product and using linearity of expectation, involving all lower-order product moments.

- **General case:** For joint moments of any number of variables, one simply substitutes the corresponding parameter list into the norm to obtain the raw product moment. By appropriate combination, central moments and cumulants of any order can be derived.

Thus, the MAO norm encapsulates all moment information of the distribution in a concise and profound manner, laying the foundation for asymptotic analysis.

### 4. Relation between Individual Marginal Distribution and Binomial Distribution

Consider a fixed individual, say individual 1. Define the indicator variables $Y_{1j} = \mathbf{1}_{1 \in S_j}$ for $j = 1, \ldots, T$. Since the subsets are drawn independently and each has size $m$, the probability that individual 1 appears in a particular subset is $p = m/N$, and events for different subsets are independent. Therefore, the number of subsets to which individual 1 belongs,

$$K_1 = \sum_{j=1}^{T} Y_{1j},$$

follows a binomial distribution:

$$\mathbb{P}(K_1 = t) = \binom{T}{t} p^t (1-p)^{T-t} \equiv \pi_t, t = 0, 1, \ldots, T.$$

This relation holds exactly for any finite $N$, without any limit. It shows that the marginal distribution of a single individual is completely described by the binomial distribution.

By linearity of expectation, we immediately obtain

$$\mathbb{E}[X_{=t}] = N\pi_t.$$

On the other hand, from the first-moment formula in the MAO norm, we have $\mathbb{E}[X_{=t}] = \| t \|_T$. Hence we obtain the identity

$$\| t \|_T = N\pi_t. \quad (3)$$

This can be understood from two perspectives: the left-hand side is a direct computation via the MAO norm, while the right-hand side is a simple expectation based on individual independence. Verifying (3) serves as a consistency check for the MAO theory (see[3]).

## 5. Covariance Between Individuals and Its Asymptotic Vanishing

Although the marginal distribution of a single individual is simple, different individuals are correlated due to sampling without replacement. For two distinct individuals $i$ and $j$, define the indicators $I_i = \mathbf{1}_{\{K_i = t\}}$ and $I_j = \mathbf{1}_{\{K_j = t\}}$. The covariance $\mathrm{Cov}(I_i, I_j)$ quantifies this dependence.

Two complementary approaches are employed to study this covariance. The first, based on the MAO inequality (Theorem 2), yields a simple order bound (Appendix proof 1):

$$\mathrm{Cov}(I_i, I_j) = O(1/N),$$

which already shows that the dependence decays rapidly with the population size.

The second approach provides an exact asymptotic expansion via a detailed analysis of the joint probability $P_t = \mathbb{P}(I_i = 1, I_j = 1)$, leading to the precise leading term:

$$\mathrm{Cov}(I_i, I_j) = \pi_t^2 \cdot \frac{1}{N} \cdot \frac{2t - Tp - \frac{t^2}{Tp}}{1 - p} + O\left(\frac{1}{N^2}\right) \quad (5.1)$$

This explicit formula reveals not only the $O(1/N)$ decay but also the precise dependence on $t$ and $p$. The detailed derivations of both the order bound and the exact expansion are provided in the Supplementary Proof 1.

## 6. Central Limit Theorem and Normal Approximation

### 6.1 Weak dependence between individuals

In the MAO model, for any two distinct individuals $i$ and $j$, let $I_i = \mathbf{1}_{\{K_i = t\}}$,

where $K_i$ is the number of subsets containing individual $i$. From the derivation in Section 5, the covariance satisfies

$$\text{Cov}(I_i, I_j) = O(1/N).$$

This implies that the correlation between any two individuals is of order $O(1/N)$.

To obtain the variance of the total count $X_{=t} = \sum_{i=1}^{N} I_i$, we use the variance decomposition:

$$\text{Var}(X_{=t}) = \sum_{i=1}^{N} \text{Var}(I_i) + \sum_{i \neq j} \text{Cov}(I_i, I_j).$$

By symmetry, all individuals have the same variance. Since each $I_i$ is a Bernoulli random variable with success probability $\pi_t$, its exact variance is $\text{Var}(I_i) = \pi_t(1 - \pi_t)$.

Moreover, all pairwise covariances are equal to $\gamma_N := \text{Cov}(I_i, I_j)$. Thus,

$$\text{Var}(X_{=t}) = N\pi_t(1 - \pi_t) + N(N-1)\gamma_N.$$

Using $\gamma_N = O(1/N)$, we obtain

$$\text{Var}(X_{=t}) = N\pi_t(1 - \pi_t) + O(N).$$

Hence, the variance grows linearly with $N$, and the leading term is $N\pi_t(1 - \pi_t)$, though a correction of order $N$ (which can be made explicit via the exact formula (5.1)) is present.

### 6.2 Applicability of the Central Limit Theorem

**Theorem 6.1 (Central Limit Theorem for MAO distribution)**

Let $N \to \infty$, $p = m/N \to p_0 \in (0,1)$, and assume $T, t$ are fixed. If $\sigma_t^2(p_0) := \pi_t(1 - \pi_t) > 0$, then

$$\frac{X_{=t} - N\pi_t}{\sqrt{\text{Var}(X_{=t})}} \Rightarrow \mathcal{N}(0,1). \qquad (6.6)$$

**Proof.**

From the results of the previous section,

$$\text{Var}(S_N) = N\pi_t(1 - \pi_t) + O(N),$$

we have $\text{Var}(S_N) \to \infty$ for sufficiently large $N$.

Let $S_N = X_{=t} = \sum_{i=1}^{N} I_i$, with $I_i = \mathbf{1}_{\{K_i = t\}}$. From Section 5, the covariance between any two distinct individuals satisfies $\text{Cov}(I_i, I_j) = O(1/N)$, and the variables are exchangeable. Consequently, we have:

- **Variance growth**: $\text{Var}(S_N) = N\pi_t(1 - \pi_t) + O(N) \to \infty$ (since $\pi_t(1 - \pi_t) > 0$);

- **Boundedness**: Each $I_i$ is Bernoulli distributed with parameter $\pi_t$, hence $|I_i| \leq 1$; therefore all moments exist and are uniformly bounded;

- **Finite-dimensional convergence**: For any fixed $k$, the vector $(I_1, \dots, I_k)$ converges in distribution to an independent and identically distributed (i.i.d.) sequence, with each component following $\text{Bernoulli}(\pi_t)$. This follows from the fact that covariances tend to zero and higher-order moments also converge to the product form (see Appendix proof 2).

These conditions are precisely those required for a central limit theorem for exchangeable sequences (see, e.g., Billingsley 1995, or Durrett 2019)[4, 5]. The theorem asserts that if $\{I_{N,i}\}$ is exchangeable, the variance tends to infinity, and every finite-dimensional distribution converges to an i.i.d. law, then the standardized sum converges in distribution to a standard normal. Hence,

$$\frac{S_N - \mathbb{E}[S_N]}{\sqrt{\text{Var}(S_N)}} \Rightarrow \mathcal{N}(0,1).$$

Substituting $\mathbb{E}[S_N] = N\pi_t$ yields (6.6).

### 7. Poisson Approximation in the Rare-Event Regime

This section considers the situation where the expectation $\lambda_N = \mathbb{E}[X_{=t}] = N\pi_t$ remains bounded (i.e., $\lambda_N = O(1)$). Under this sparse condition, the event that an individual belongs to exactly $t$ subsets is rare, and the dependence among individuals becomes negligible; consequently, the limiting distribution of $X_{=t}$ is Poisson.

#### 7.1 Statement of the Poisson Limit

Let $p = m/N$ vary with $N$ such that $\lambda_N := N\pi_t \to \lambda \in [0, \infty)$, where $\pi_t = \binom{T}{t} p^t (1-p)^{T-t}$. Define

$$\Delta_N := \sum_{1 \leq i \neq j \leq N} |\mathbb{P}(I_i = 1, I_j = 1) - \pi_t^2|.$$

**Theorem 7.1 Poisson Limit Theorem.**

Let $\lambda_N = \mathbb{E}[X_{=t}] \to \lambda \in [0, \infty)$, which implies $\mathbb{E}[X_{=t}] \to 0$ when $N \to \infty$. Then

$$X_{=t} \Rightarrow \text{Poisson}(\lambda). \qquad (7.1)$$

*Proof.* Applying the Chen–Stein Poisson approximation method[6, 7], the total variation distance satisfies

$$d_{\text{TV}}(\mathcal{L}(X_{=t}), \text{Poisson}(\lambda_N)) \leq \sum_{i=1}^{N} \pi_t^2 + \sum_{i \neq j} |\mathbb{E}[I_i I_j] - \mathbb{E}[I_i]\mathbb{E}[I_j]| = N\pi_t^2 + \Delta_N. \qquad (7.2)$$

From $\lambda_N = N\pi_t \to \lambda$ we have $N\pi_t^2 = \lambda_N^2/N \to 0$. Together with $\Delta_N \to 0$, the right-hand side tends to zero, so $d_{\text{TV}}(\mathcal{L}(X_{=t}), \text{Poisson}(\lambda_N)) \to 0$. Since $\lambda_N \to \lambda$, we obtain (7.1).

### 7.2 Verification of the Second-Order Condition via the MAO Norm

To apply Theorem 7.1, we need to verify that $\Delta_N \to 0$ when $\lambda_N = N\pi_t$ is bounded (i.e., $\lambda_N = O(1)$). The key observation is that in the sparse regime where $\lambda_N \to 0$ (or remains bounded), the MAO distribution exhibits a fundamental property: the variance of $X_{=t}$ is asymptotically equal to its mean, and the difference is of higher order relative to the mean. More precisely, from the MAO theory it is known that

$$\text{Var}(X_{=t}) = \lambda_N - \delta_N, \text{with } \delta_N = o(\lambda_N) \text{ as } \lambda_N \to 0,$$

and moreover $\text{Var}(X_{=t}) < \lambda_N$ (the MAO inequality). This property follows from the exact moment formulas and has been established in previous work[3]. We now use it to bound $\Delta_N$.

Recall the definition of $\Delta_N$:

$$\Delta_N = \sum_{1 \leq i \neq j \leq N} |\mathbb{P}(I_i = 1, I_j = 1) - \pi_t^2|.$$

By exchangeability, all off-diagonal terms are equal, so $\Delta_N = N(N-1) |P_t - \pi_t^2|$ with $P_t = \mathbb{P}(I_1 = 1, I_2 = 1)$. The variance of $X_{=t}$ can be expressed as

$$\text{Var}(X_{=t}) = \sum_{i=1}^{N} \text{Var}(I_i) + \sum_{i \neq j} \text{Cov}(I_i, I_j) = N\pi_t(1 - \pi_t) + N(N-1)(P_t - \pi_t^2).$$

Solving for $P_t - \pi_t^2$ gives

$$P_t - \pi_t^2 = \frac{\text{Var}(X_{=t}) - N\pi_t(1-\pi_t)}{N(N-1)}.$$

Now, using $N\pi_t = \lambda_N$ and $N\pi_t^2 = \lambda_N^2/N$, we have

$$N\pi_t(1 - \pi_t) = \lambda_N - \frac{\lambda_N^2}{N}.$$

Hence

$$P_t - \pi_t^2 = \frac{\text{Var}(X_{=t}) - \lambda_N + \frac{\lambda_N^2}{N}}{N(N-1)}.$$

Substituting $\text{Var}(X_{=t}) = \lambda_N - \delta_N$ (with $\delta_N = o(\lambda_N)$) yields

$$P_t - \pi_t^2 = \frac{-\delta_N + \frac{\lambda_N^2}{N}}{N(N-1)}.$$

Therefore,

$$\Delta_N = N(N-1) \mid P_t - \pi_t^2 \mid = \mid -\delta_N + \frac{\lambda_N^2}{N} \mid \leq \delta_N + \frac{\lambda_N^2}{N}.$$

Since $\lambda_N = O(1)$, we have $\lambda_N^2/N \to 0$. Moreover, $\delta_N = o(\lambda_N)$ implies $\delta_N \to 0$ (because $\lambda_N$ is bounded). Consequently, $\Delta_N \to 0$ as $N \to \infty$.

Thus, the second-order condition $\Delta_N \to 0$ is automatically satisfied in the sparse regime where $\lambda_N$ remains bounded. This completes the verification of the hypothesis of Theorem 7.1, and together with $\lambda_N \to \lambda$ yields the Poisson limit

$$X_{=t} \Rightarrow \text{Poisson}(\lambda).$$

**Remark.** The above argument relies only on the asymptotic equality of mean and variance (and the fact that the variance is smaller than the mean) which is a known property of the MAO distribution when the expected count is small. This elegant approach avoids explicit expansion of the MAO norm and highlights the essential role of the moment relationship in the Poisson approximation.

### 7.3 Parameter Scales and Practical Interpretation

Because $\lambda_N = \mathbb{E}[X_{=t}] = N\pi_t$ completely determines the applicability of the Poisson approximation, we can classify the limiting behaviour directly by the asymptotic behaviour of the mean:

- If $\mathbb{E}[X_{=t}] \to \lambda \in (0, \infty)$, and $N \to \infty$, then $X_{=t}$ converges to a non-degenerate Poisson distribution $\text{Poisson}(\lambda)$.

- If $\mathbb{E}[X_{=t}] \to 0$ (which can occur, for example, when $N$ is fixed and $T \to \infty$, or when $m$ grows sufficiently slowly), then $X_{=t}$ converges in distribution to a degenerate Poisson distribution (i.e., a point mass at zero).

- If $\mathbb{E}[X_{=t}] \to N$ (i.e., the expectation approaches the total population size), then $N - X_{=t}$ converges to a Poisson distribution (degenerate at zero), and $X_{=t}$ is approximately $N$ minus a Poisson variable. This is the mirror image of the case $\mathbb{E}[X_{=t}] \to 0$.

- In the intermediate regime where $\mathbb{E}[X_{=t}]$ is large but bounded away from both $0$ and $N$ (i.e., $\mathbb{E}[X_{=t}]/N \to c \in (0,1)$), the normal approximation (Section 6) applies.

This simple criterion avoids complicated parameter expansions and highlights the central role of the expectation in limit theory. For example, when $p \to 0$ we have $\lambda_N \sim \binom{T}{t} m^t / N^{t-1}$, but one need not memorize this formula; it suffices to check whether $\lambda_N$ tends to a constant, to zero, to $N$, or to an intermediate value.

### 7.4 Summary

The essence of the Poisson approximation lies in the superposition of rare events: when each individual's probability $\pi_t$ of being a "success" is small, and the expected count $\lambda_N$ is moderate, the total count $X_{=t}$ behaves like a sum of independent Bernoulli variables, even though weak dependence exists; its effect becomes negligible in the limit. This part has rigorously established this fact for the MAO distribution using its intrinsic moment relations, and has reduced the condition for the Poisson approximation to the boundedness of the mean, providing clear guidance for practical applications.

## 8. Numerical Verification

To visually assess the goodness of fit between the MAO distribution and its normal and Poisson approximations, we conduct comparisons under various parameter regimes using both exact probability calculations and Monte Carlo simulations. The focus is on the count variables $X_{=t}$ and their cumulative counterparts $X_{\geq t} = \sum_{s=t}^{T} X_{=s}$. Throughout this section, the number of subsets is fixed at $T = 5$, and the proportion $p = m/N$ is chosen according to the scenario. To provide a comprehensive view, each figure is organized into two large panels: the left panel compares the MAO distribution with the normal approximation, and the right panel compares it with the Poisson approximation. Within each large panel, two columns display $X_{=t}$ (left column) and $X_{\geq t}$ (right column) for $t = 2,3,4,5$ from top to bottom.

### 8.1 Exact Distribution Comparison for Small Samples (Equal Subset Sizes)

We first consider a small population $N = 100$ with equal subset sizes $m = 20$ (i.e., $p = 0.2$). Moments computed via the MAO norm formulas, an exact recursive algorithm, and Monte Carlo simulations ($10^5$ runs) agree perfectly, validating the MAO moment formulas (see Table S1 in the supplement).

Figure 1 compares the exact probability distributions with the normal and Poisson approximations. The left large panel shows the normal comparison, the right large panel the Poisson comparison; within each, the left column displays $X_{=t}$ and the right column $X_{\geq t}$ for $t = 2,3,4,5$ (top to bottom). We observed that when the mean is large (e.g., $t = 2,3$), the normal curve closely matches the exact distribution. As the mean decreases (increasing $t$), the normal fit deteriorates while the Poisson approximation becomes increasingly accurate; for very small means (e.g., $t = 5$, mean $\approx 0.032$), the Poisson is nearly perfect. The cumulative variables $X_{\geq t}$ exhibit the same pattern.

### 8.2 Simulated Distribution Comparison for Large Samples (Equal Subset Sizes)

Now take $N = 5000$ and $m = 1000$ (still $p = 0.2$). Exact computation is infeasible, so we perform $10^5$ Monte Carlo simulations. The empirical moments agree closely with the theoretical ones (Table S2).

Figure 2 adopts the same layout as Figure 1. The left panel shows that for all $t = 2,3,4,5$, the normal density aligns almost perfectly with the histograms, confirming the central limit theorem. In the right panel, the Poisson approximation is excellent when the mean is small (e.g., $t = 4,5$, with means $\approx 32$ and $1.6$), while for larger means ($t = 2,3$, means $\approx 1024$ and $256$) the Poisson is omitted as it visually coincides with the normal and is less informative. The cumulative variables $X_{\geq t}$ exhibit identical behavior. These observations mirror the pattern seen in the exact small-sample case (Figure 1): normal works well for moderately large means, while Poisson excels in the sparse regime where the mean is small.

### 8.3 Unequal Subset Sizes (Exact Distribution)

In practice subset sizes may vary. Take $N = 100$ with **m** = [5,20,40,70,30] (range from 5 to 70) and compute the exact distribution via a recursive algorithm (moments in Table S3). Figure 3 (same layout) demonstrates that, despite the substantial variation in subset sizes, the same pattern emerges: the normal approximation works well when the mean is large (e.g., $t = 2,3$), and the Poisson approximation is excellent when the mean is small (e.g., $t = 4,5$). This indicates that the asymptotic behavior is robust to the equality assumption and is primarily governed by the magnitude of the mean.

### 8.4 summary

The numerical experiments perfectly validate the MAO norm moment formulas: in Table S1, the formula and exact values coincide to machine precision, confirming the theoretical soundness of the MAO norm as a moment-generating tool. The simulation results, with only minor sampling fluctuations, further corroborate the structural integrity of the MAO distribution.

Remarkably, even when subset sizes are highly unequal (e.g., ranging from 5 to 70), the distributions of $X_{=t}$ still conform closely to either a normal or a Poisson approximation, with the choice governed solely by the magnitude of the mean: large mean favors normal, small mean favors Poisson. This suggests that the inherent symmetry or exchangeability of the MAO design may play a key role. Although the individual correlations depend on the specific $m_i$ values, the aggregate behavior of the count variables is dominated by their means. The exchangeability among individuals (despite differing subset sizes) likely averages out higher-order fluctuations, rendering the influence of unequal $m_i$ negligible. A deeper explanation lies in the fact that, for fixed $T$ and sufficiently large $N$, $X_{=t}$ is a sum of many weakly dependent random variables; the central limit theorem thus overwhelms any local perturbations due to varying subset sizes. This insight provides practical guidance: in applications, one may safely use normal or Poisson approximations for MAO counts as long as the mean is appropriate, regardless of whether the subset sizes are equal.

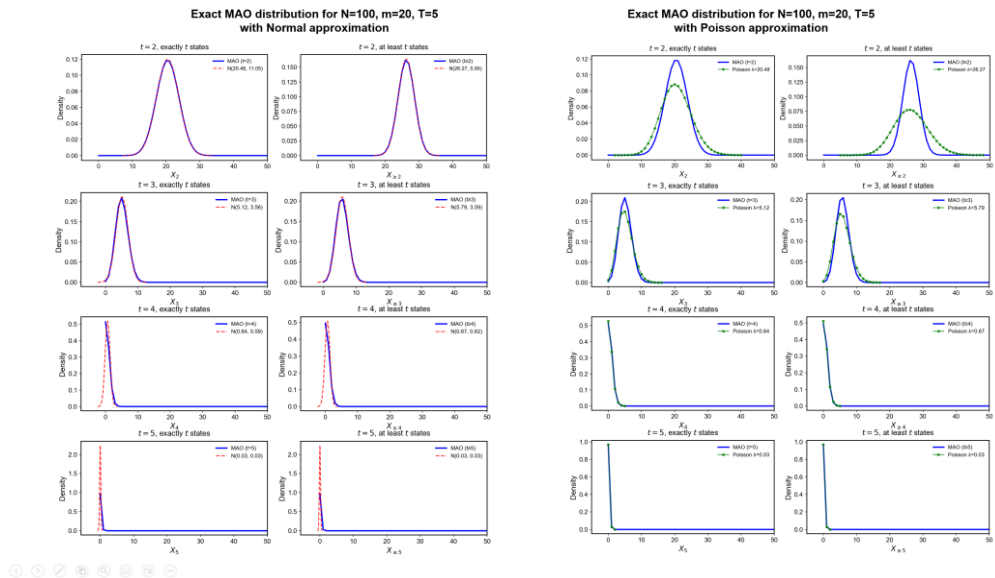

**Figure 1: Exact distributions vs. approximations for small samples with equal subset sizes ($N = 100, m = 20, T = 5$).**

Left large panel: comparison with the normal approximation. Subplots in the first column (top to bottom: $t = 2,3,4,5$) show the exact probabilities (blue vertical bars) of $X_{=t}$ together with the normal density (red solid curve) using the exact moments from Table S1 (e.g., $\mu_2 \approx 20.48, \sigma_2^2 \approx 11.05$; $\mu_3 \approx 5.12, \sigma_3^2 \approx 3.56$; $\mu_4 \approx 0.64, \sigma_4^2 \approx 0.59$; $\mu_5 \approx 0.032, \sigma_5^2 \approx 0.032$). Subplots in the second column show the exact probabilities of $X_{\geq t}$ with the normal density. Right large panel 3rd and 4th columns:

comparison with the Poisson approximation. Subplots in the 3rd column show the exact probabilities of $X_{=t}$ together with the Poisson mass (green dots) with parameter $\lambda_t = \mu_t$; subplots in the 4th column show the exact probabilities of $X_{\geq t}$ with the Poisson mass. The plots illustrate that normal fits are excellent when the mean is large, while Poisson fits are superior for rare events.

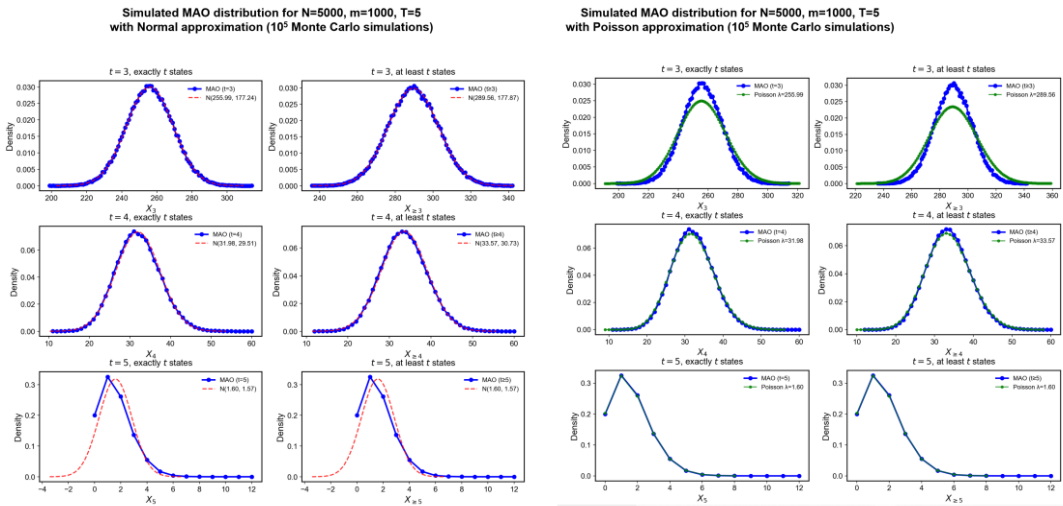

**Figure 2: Simulated distributions vs. approximations for large samples with equal subset sizes ($N = 5000, m = 1000, T = 5$).**

Based on $10^5$ Monte Carlo simulations, the layout mirrors Figure 1. Left large panel: histograms (grey bars) of $X_{=t}$ (left column) and $X_{\geq t}$ (right column) overlaid with normal density (red). Theoretical moments from Table S2 (e.g., $\mu_2 \approx 1024, \sigma_2^2 \approx 552.15$; $\mu_3 \approx 256, \sigma_3^2 \approx 177.37$; $\mu_4 \approx 32, \sigma_4^2 \approx 29.49$; $\mu_5 \approx 1.6, \sigma_5^2 \approx 1.59$). Right large panel: histograms overlaid with Poisson mass (green dots).

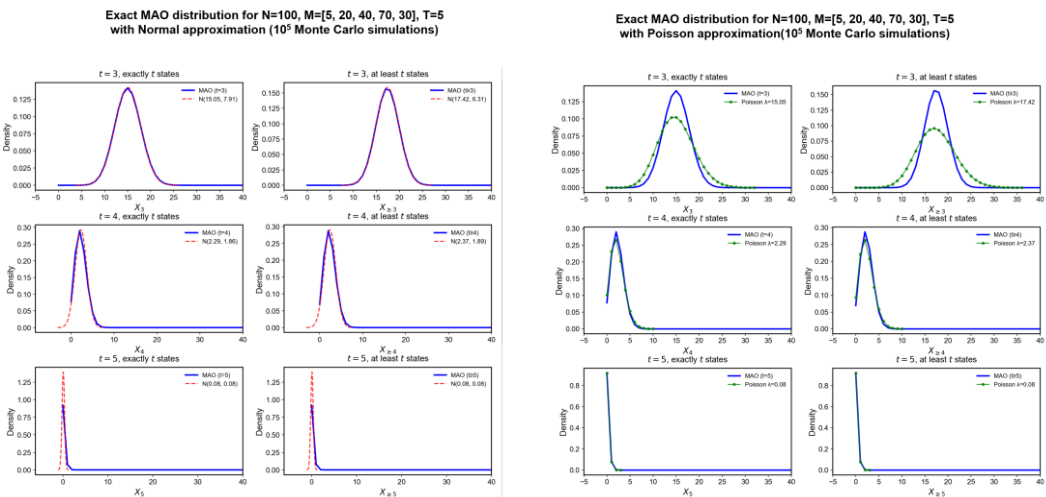

**Figure 3: Exact distributions vs. approximations for unequal subset sizes ($N =$

$100, \mathbf{m} = [5,20,40,70,30], T = 5$).

Layout identical to Figure 1, with exact distributions obtained via a recursive algorithm (moments in Table S3: $\mu_2 \approx 37.27, \sigma_2^2 \approx 20.98$; $\mu_3 \approx 15.05, \sigma_3^2 \approx 7.91$; $\mu_4 \approx 2.29, \sigma_4^2 \approx 1.86$; $\mu_5 \approx 0.084, \sigma_5^2 \approx 0.082$). Despite the substantial variation in subset sizes, the same pattern emerges: normal approximation works well when the mean is large (e.g., $t = 2,3$), and Poisson approximation is excellent when the mean is small (e.g., $t = 4,5$). This demonstrates the robustness of the asymptotic behavior.

## 9. Conclusion and Discussion

This paper has provided a comprehensive analysis of the Multi-set Allocation Occupancy (MAO) distribution—the joint distribution of the count vector $\mathbf{X} = (X_{=0}, X_{=1}, \ldots, X_{=T})$ arising from $T$ independent rounds of sampling without replacement, each of size $m$, from a population of $N$ individuals. By introducing the **MAO norm** as a unifying algebraic tool, we have obtained closed-form expressions for all factorial moments, rigorously established the asymptotic behavior of $X_{=t}$ under two canonical regimes, and elucidated the connection between the MAO distribution and the classical binomial, Poisson and normal distributions. The main theoretical findings are summarised in Sections 9.1–9.7. Building on these results, we now reflect on the broader significance of the MAO distribution as a null model for detecting non-random aggregation in multi-state data, and outline its potential across several scientific disciplines.

### 9.1 Marginal Distribution and Global Coupling

For any fixed individual $i$, the number of times it is selected in the $T$ rounds, $K_i = \sum_{k=1}^{T} \mathbf{1}_{i \in S_k}$, follows an exact binomial distribution $\text{Bin}(T, p)$ with $p = m/N$. Thus, the novelty of the MAO distribution does not lie in the marginal behavior but in the global coupling (negative correlation) among individuals induced by the fixed sample size per round. This coupling causes the variance of the count variable $X_{=t}$ to deviate from that of an independent binomial model by a term of order $O(1)$.

### 9.2 Exact Order of Weak Dependence

For any two distinct individuals, define the indicator $I_i = \mathbf{1}_{\{K_i = t\}}$. Then the covariance satisfies

$$\text{Cov}(I_i, I_j) = O(1/N).$$

This decay rate means that, as $N \to \infty$, any finite collection of individuals becomes asymptotically independent; however, for the total count $X_{=t} = \sum_{i=1}^{N} I_i$, the cumulative contribution $N(N-1)\text{Cov}(I_i, I_j)$ is of order $O(N)$, which is of the same order as the leading variance term and therefore cannot be ignored. This explains why the variance of $X_{=t}$ must include a second-order correction and why the limiting distribution cannot be simply borrowed from an independent model.

### 9.3 Normal Approximation: Macroscopic Fluctuations under Fixed $p$

When $N \to \infty$ with $p = m/N$ fixed in $(0, 1)$ and $T, t$ fixed, $X_{=t}$ obeys a central limit theorem:

$$\frac{X_{=t} - N\pi_t}{\sqrt{\text{Var}(X_{=t})}} \Rightarrow \mathcal{N}(0, 1),$$

Considering that $N\pi_t = E(X_{=t})$, under the "dense sampling" regime, the fluctuations of $X_{=t}$ scale as $\sqrt{N}$ and can be approximated by a normal distribution:

$$X_{=t} \approx \mathcal{N}(E(X_{=t}), \text{Var}(X_{=t})).$$

### 9.4 Poisson Approximation: Rare-Event Regime

The Poisson approximation applies when the expectation $\lambda_N = \mathbb{E}[X_{=t}]$ remains bounded as $N \to \infty$. More precisely, if $\lambda_N \to \lambda \in [0, \infty)$, then $X_{=t}$ converges in distribution to a Poisson distribution with parameter $\lambda$ (with the convention that $\text{Poisson}(0)$ is the point mass at zero). In the MAO model, this is guaranteed by the decay of covariances established in Section 7, which ensures the weak dependence condition required for Poisson approximation (e.g., via the Chen–Stein method).

### 9.5 Unifying Role of the MAO Norm

The technical thread of this paper centers on the MAO norm: it provides closed-form expressions for all factorial moments, systematizes moment calculations, and enables the derivation of limit distributions via moment convergence under different parameter regimes. The MAO norm is not merely a computational tool but a bridge connecting the exact combinatorial structure to asymptotic theory.

### 9.6 Practical Guidelines

For finite $N$, the choice of approximation for $X_{=t}$ can be guided by the magnitude of its mean $\mu = \mathbb{E}[X_{=t}]$ relative to the population size $N$:

1. **When $\mu$ is very small** (compared to $N$), the distribution of $X_{=t}$ is well approximated by a Poisson distribution with parameter $\mu$.

2. **When $\mu$ is very close to $N$**, the complement $N - X_{=t}$ is approximately

Poisson, so $X_{=t}$ behaves like $N$ minus a Poisson variable.

3. **When $\mu$ is neither too small nor too close to $N$** (i.e., $\mu/N$ is moderate), the normal approximation $\mathcal{N}(\mu, \text{Var}(X_{=t}))$ is appropriate.

These simple criteria, based solely on the size of the mean, avoid complicated parameter expansions and reflect the essential role of the expectation in limit theory. For instance, when $p \to 0$ we have $\mu \sim \binom{T}{t} m^t / N^{t-1}$; one need not memorize this formula—it suffices to check whether $\mu$ is tiny, intermediate, or close to $N$.

### 9.7 The MAO Distribution as a Universal Null Model for Aggregation

At its core, the MAO distribution describes the random allocation of individuals to multiple states (the subsets) under the simplest possible mechanism: each state is formed by an independent, uniform, without-replacement sample of fixed size. No further structure – such as interactions between individuals, preferences for certain states, or heterogeneity among states – is assumed. Consequently, the resulting distribution of $X_{=t}$ (the number of individuals that occupy exactly $t$ states) serves as a **rigorous null model** against which empirical data can be compared. If an observed count deviates significantly from the MAO prediction, it signals the presence of non-random features, such as:

- **Aggregation**: individuals tend to co-occur in the same states more often than expected by chance;

- **Repulsion**: individuals avoid co-occurrence, leading to fewer multi-state individuals than expected.

The need for such a null model is especially acute because the simple binomial approximation (treating individuals as independent) ignores the weak negative correlations induced by sampling without replacement; these correlations, although of order $O(1/N)$, accumulate and affect the variance at the $O(N)$ level, as shown in Section 9.2. The MAO distribution captures this subtle dependence exactly, making it the correct benchmark for any finite population.

### 9.8 Cross-Disciplinary Perspectives

The abstract structure of the MAO model – a finite set of individuals, a collection of states (subsets), and the count of individuals with a given multiplicity – appears naturally in numerous fields. By interpreting the mathematical objects appropriately, one can translate the MAO results into concrete scientific insights.

#### 9.8.1 Ecology: Species Co-occurrence Patterns

In ecology, let the individuals be species and the states be sampling sites (e.g., quadrats, traps). Each site records the presence of a fixed number $m$ of species (a typical situation when sampling effort is standardized). Then $X_{=t}$ is the number of species that occur in exactly $t$ sites. Ecologists have long debated whether observed species co-occurrences reflect random assembly or structured communities (e.g., competition, mutualism)[8, 9]. The MAO distribution provides a mathematically sound null hypothesis: under random placement, the expected number of species with $t$ occurrences is $N\pi_t$, and the full distribution is known (exactly for small $N$, approximately normal or Poisson for larger $N$). A significant excess of species with high $t$ may indicate positive associations (e.g., shared habitat preferences), while a deficit may suggest negative interactions. Our moment formulas and asymptotic approximations enable rigorous testing, even when the number of sites is modest.

### 9.8.2 Sociology and Network Science: Multiple Group Membership

Consider a population of individuals and a collection of social groups (clubs, online communities, project teams). If each group is formed by randomly selecting $m$ members from the population (without replacement), then $X_{=t}$ counts how many individuals belong to exactly $t$ groups. In reality, group membership is often non-random: some individuals are "joiners" who belong to many groups, while others remain isolated[10, 11]. The MAO distribution offers a baseline against which such social activity can be assessed. For example, a significantly larger than expected number of individuals with $t \geq 2$ might indicate the presence of influential hubs or overlapping social circles. Conversely, a deficit could signal fragmentation. The same logic applies to overlapping communities in networks: a node's membership in multiple communities can be tested against the MAO null model.

### 9.8.3 Physics: Multiple Occupancy of Energy Levels

In statistical mechanics, a system of $N$ particles distributed over $T$ energy levels, each capable of holding up to $m$ particles, is a classic setting (e.g., in the context of Bose–Einstein condensation or Pauli exclusion)[12-14]. Here $X_{=t}$ is the number of levels occupied by exactly $t$ particles. If the particles are non-interacting and the levels are chosen randomly with a fixed capacity (a "canonical" type of sampling), the MAO distribution emerges as the exact occupancy distribution. Deviations from this random-occupancy prediction would signal interactions (attraction or repulsion) between particles. Our results provide the necessary statistical tools to quantify such deviations, especially in the dilute limit where the Poisson approximation applies, or in the dense limit where normal fluctuations dominate.

### 9.8.4 Other Domains

- **Bioinformatics**: In functional genomics, genes (individuals) may be annotated to multiple pathways (states). The number of genes shared by $t$ pathways can be modelled by MAO; enrichment of genes in many pathways might indicate pleiotropy[15, 16].

- **Epidemiology**: Individuals may be co-infected by several pathogens (states)[17, 18]. The count of individuals with $t$ infections can be compared with the MAO null to detect clustering of infections.

- **Library Science**: Documents (individuals) may be assigned to multiple subject categories (states)[19, 20]; the distribution of documents across categories can be benchmarked.

### 9.9 From Theory to Practice: A Roadmap

The theoretical developments of this paper translate directly into practical recommendations for researchers in any of the above fields:

- **Estimate parameters**: From the data, obtain $N$ (total individuals), $T$ (number of states), and the average state size $m$ (or the individual probabilities $p = m/N$).

- **Compute expected counts**: Use $\pi_t = \binom{T}{t} p^t (1-p)^{T-t}$ to obtain the expected number of individuals with exactly $t$ states.

- **Assess variability**: For small $N$, use exact MAO probabilities (available via recursive algorithms or pre-computed tables). For large $N$, employ the normal approximation when $N\pi_t$ is large (say >10) and the Poisson approximation when $N\pi_t$ is of order 1. The variance formulas derived from the MAO norm guarantee that the approximations are correctly scaled.

- **Test hypotheses**: Construct confidence intervals or p-values based on the chosen approximation. A significant deviation points to non-random structure, warranting further investigation.

### 9.10 Limitations and Future Work

While the MAO framework is remarkably versatile, several extensions would increase its applicability:

- **Variable subset sizes**: This paper assumed equal subset sizes $m$. In many applications, the sizes may differ. Preliminary numerical experiments (Section 9.3) suggest that the asymptotic behaviour is robust to moderate heterogeneity, but a rigorous theory for unequal $m_i$ remains to be developed.

- **Growing $T$ or $t$**: When the number of states $T$ grows with $N$ (e.g., in high-throughput genomics), new phase transitions may appear. The current results, valid for fixed $T$, need to be extended to the regime where $T = T_N \to \infty$.

- **Joint distribution of $(X_{=0}, \ldots, X_{=T})$**: For many applications, the full vector is of interest. Its limit (a degenerate multivariate normal or a multivariate Poisson) could be derived using similar moment-based techniques, providing tools for simultaneous inference across all $t$.

- **Sharp error bounds**: The approximations presented are first-order; deriving explicit bounds on the total variation or Kolmogorov distance (e.g., via Stein's method) would give practitioners precise guarantees for finite $N$.

### 9.11 Concluding Remarks

The MAO distribution occupies a unique position at the intersection of probability theory, combinatorics, and applied statistics. It arises from one of the simplest and most ubiquitous sampling designs – independent without-replacement draws – yet its exact distribution is far from trivial. By harnessing the power of the MAO norm, we have unlocked its moment structure and asymptotic behavior, revealing a beautiful connection with the three pillars of classical probability: the binomial, Poisson and normal laws. More importantly, the MAO distribution provides a principled null model for detecting aggregation in any system where individuals can possess multiple states. We hope that the theoretical framework established here will stimulate both further mathematical research and widespread application across the natural, social and computational sciences.

**Appendix Proof 1: Covariance Decay via the MAO Inequality:** $\mathrm{Cov}(I_i, I_j) = O(1/N)$

The MAO inequality (Theorem 2) states that for any list of parameters $(p_1, \ldots, p_r)$ when $p_i$ are close enough,

$$\prod_{j=1}^{r} \| p_j \|_T \geq \| (p_1, \ldots, p_r) \|_T.$$

Taking $r = 2$ and $p_1 = p_2 = t$ yields

$$\| t \|_T^2 \geq \| (t, t) \|_T.$$

Since $\| t \|_T = N\pi_t$, we obtain

$$\| (t, t) \|_T \leq N^2 \pi_t^2.$$

Recall that for two distinct individuals,

$$P_t := \mathbb{P}(I_1 = 1, I_2 = 1) = \frac{\| (t, t) \|_T}{N(N-1)}.$$

Thus,

$$P_t \leq \frac{N^2 \pi_t^2}{N(N-1)} = \pi_t^2 \cdot \frac{N}{N-1} = \pi_t^2 + \frac{\pi_t^2}{N-1}.$$

Hence,

$$P_t - \pi_t^2 \leq \frac{\pi_t^2}{N-1} = O\left(\frac{1}{N}\right). \qquad (1)$$

To obtain a lower bound, we use the variance decomposition:

$$\mathrm{Var}(X_{=t}) = N\pi_t(1 - \pi_t) + N(N-1)(P_t - \pi_t^2).$$

Since variance is non-negative,

$$N(N-1)(P_t - \pi_t^2) \geq -N\pi_t(1 - \pi_t),$$

which gives

$$P_t - \pi_t^2 \geq -\frac{\pi_t(1 - \pi_t)}{N-1} = O\left(\frac{1}{N}\right). \qquad (2)$$

Combining (1) and (2), we conclude

$$| P_t - \pi_t^2 | = O\left(\frac{1}{N}\right),$$

and therefore

$$\text{Cov}(I_i, I_j) = P_t - \pi_t^2 = O\left(\frac{1}{N}\right).$$

This simple argument, relying only on the MAO inequality and basic variance properties, establishes the $O(1/N)$ decay of the covariance without requiring an exact asymptotic expansion. For the precise coefficient in front of $1/N$, however, one must resort to the full expansion of the MAO norm as derived earlier.

**Appendix Proof 2:**

**Proof of Asymptotic Independence and Identical Distribution of $(I_1, \ldots, I_k)$**

Consider a population of size $N$ and $T$ subsets $S_1, \ldots, S_T$, each of size $m$, drawn independently without replacement uniformly from the $N$ individuals. Let $p = m/N$ and assume that $p$ is fixed and satisfies $0 < p < 1$ (i.e., there exists $\epsilon > 0$ such that $\epsilon \le p \le 1 - \epsilon$, ensuring that coefficients appearing in expansions are bounded). For a fixed integer $t$ ($0 \le t \le T$), define the indicator variables $I_i = \mathbf{1}_{\{K_i = t\}}$, where $K_i$ is the number of subsets containing individual $i$. The marginal probability is

$$\pi_t = \binom{T}{t} p^t (1-p)^{T-t}.$$

For an arbitrary fixed positive integer $k$, set

$$P_k = \mathbb{P}(I_1 = 1, \ldots, I_k = 1).$$

**1. Relation between the MAO Norm and Moments**

From the MAO theory (see Section 3 of this paper), for the parameter list $(t, \ldots, t)$ (repeated $k$ times), we have

$$\mathbb{E}[(X_{=t})_k] = \| (t, \ldots, t) \|_T,$$

where $(X_{=t})_k = X_{=t}(X_{=t} - 1) \cdots (X_{=t} - k + 1)$ is the falling factorial, and the norm is defined by

$$\| (t, \ldots, t) \|_T = \frac{1}{(N)_k^{T-1}} \sum_{\substack{A_1, \ldots, A_k \subseteq [T] \\ |A_j| = t}} g(A_1, \ldots, A_k).$$

Here $(N)_k = N(N-1) \cdots (N-k+1)$, and

$$g(A_1, \ldots, A_k) = \prod_{i=1}^{T} (m)_{k_i} (N-m)_{k-k_i},$$

where $k_i = \sum_{j=1}^{k} \mathbf{1}_{i \in A_j}$ counts how many of the selected subsets contain index $i$. Each $A_j$ is a subset of $[T]$ of size $t$; thus the total number of tuples $(A_1, \ldots, A_k)$ is

$$M = \binom{T}{t}^k,$$

a constant independent of $N$. This expression directly reflects the structure of independent sampling without replacement in each of the $T$ rounds (multiple hypergeometric model).

On the other hand, by exchangeability of individuals,

$$\mathbb{E}[(X_{=t})_k] = \sum_{i_1 \neq \cdots \neq i_k} \mathbb{E}[I_{i_1} \cdots I_{i_k}] = (N)_k\, P_k,$$

because there are $(N)_k$ ordered $k$-tuples and all expectations are equal. Consequently,

$$P_k = \frac{\|(t,\ldots,t)\|_T}{(N)_k} = \frac{\sum_{A_1,\ldots,A_k} g(A_1,\ldots,A_k)}{(N)_k^T}. \tag{1}$$

## 2. Uniform Asymptotic Expansion of the MAO Function

Insert $m = pN$ into the expression for $g$. For a fixed $k$, each $k_i$ satisfies $0 \leq k_i \leq k$ and $\sum_{i=1}^T k_i = kt$. Using the expansion of the falling factorial, for any $0 \leq r \leq k$ we have

$$(pN)_r = (pN)^r \left(1 - \frac{\binom{r}{2}}{pN} + O\left(\frac{1}{N^2}\right)\right),$$

$$((1-p)N)_r = ((1-p)N)^r \left(1 - \frac{\binom{r}{2}}{(1-p)N} + O\left(\frac{1}{N^2}\right)\right)$$

where the constants in the $O(N^{-2})$ terms depend only on $k$ and on the lower bound of $p$ (since $p$ is fixed, the bound is uniform). Hence, for each fixed $i$,

$$(m)_{k_i}(N-m)_{k-k_i}$$

$$= (pN)^{k_i}((1-p)N)^{k-k_i}\left(1 - \frac{\binom{k_i}{2}}{pN} + O\left(\frac{1}{N^2}\right)\right)\left(1 - \frac{\binom{k-k_i}{2}}{(1-p)N} + O\left(\frac{1}{N^2}\right)\right)$$

Multiplying over $i = 1,\ldots,T$ and using $\prod_{i=1}^T (pN)^{k_i}((1-p)N)^{k-k_i} = p^{kt}(1-p)^{kT-kt}N^{kT}$ (because $\sum k_i = kt$), we obtain

$$g(A_1,\ldots,A_k) = p^{kt}(1-p)^{kT-kt}N^{kT}\prod_{i=1}^T\left[1 - \frac{\binom{k_i}{2}}{pN} - \frac{\binom{k-k_i}{2}}{(1-p)N} + O\left(\frac{1}{N^2}\right)\right].$$

Here we have combined the two first-order terms and absorbed all cross-terms (e.g., $\frac{\binom{k_i}{2}}{pN} \cdot \frac{\binom{k_j}{2}}{pN}$ for $i \neq j$) as well as the higher-order errors into the $O(N^{-2})$ term. Since $T$ is fixed, when expanding the product all cross-terms are $O(N^{-2})$.

Consequently,

$$g(A_1, \ldots, A_k) = p^{kt}(1-p)^{kT-kt}N^{kT}\left[1 - \frac{1}{N}\left(\sum_{i=1}^{T}\frac{\binom{k_i}{2}}{p} + \sum_{i=1}^{T}\frac{\binom{k-k_i}{2}}{1-p}\right) + O\left(\frac{1}{N^2}\right)\right], \quad (2)$$

where the $O(N^{-2})$ term is uniform over all tuples $(A_1, \ldots, A_k)$ because all quantities involved are bounded.

### 3. Summation of the Leading Term and the Remainder

Summing (2) over all $M = \binom{T}{t}^k$ tuples, we set

$$S = \sum_{A_1,\ldots,A_k} g(A_1, \ldots, A_k).$$

Then

$$S = p^{kt}(1-p)^{kT-kt}N^{kT}\left[M - \frac{1}{N}\sum_{A_1,\ldots,A_k}\left(\sum_{i=1}^{T}\frac{\binom{k_i}{2}}{p} + \sum_{i=1}^{T}\frac{\binom{k-k_i}{2}}{1-p}\right) + O\left(\frac{M}{N^2}\right)\right].$$

Since $\binom{k_i}{2} \leq \binom{k}{2}$ and $\binom{k-k_i}{2} \leq \binom{k}{2}$, and $1/p, 1/(1-p)$ are bounded (because $p$ is fixed and bounded away from $0$ and $1$), we have

$$\left|\sum_{A_1,\ldots,A_k}\left(\sum_{i=1}^{T}\frac{\binom{k_i}{2}}{p} + \sum_{i=1}^{T}\frac{\binom{k-k_i}{2}}{1-p}\right)\right| \leq M \cdot T \cdot \binom{k}{2} \cdot \left(\frac{1}{p} + \frac{1}{1-p}\right) = O(1),$$

where the constant depends on $T, t, k, p$. Therefore,

$$S = p^{kt}(1-p)^{kT-kt}N^{kT}\left[M + O\left(\frac{1}{N}\right)\right].$$

Using $\pi_t = \binom{T}{t}p^t(1-p)^{T-t}$, we have $p^{kt}(1-p)^{kT-kt} = \pi_t^k/\binom{T}{t}^k$, and recalling that $M = \binom{T}{t}^k$, we obtain

$$S = \pi_t^k N^{kT}\left[1 + O\left(\frac{1}{N}\right)\right]. \quad (3)$$

### 4. The Denominator

From (1), $P_k = S/(N)_k^T$. For a fixed $k$, the falling factorial expands as

$$(N)_k = N^k \prod_{j=0}^{k-1}\left(1-\frac{j}{N}\right) = N^k\left(1 - \frac{k(k-1)}{2N} + O\left(\frac{1}{N^2}\right)\right) = N^k\left(1+O\left(\frac{1}{N}\right)\right)$$

Thus,

$$(N)_k^T = N^{kT}\left(1+O\left(\frac{1}{N}\right)\right)$$

Inserting (3) yields

$$P_k = \frac{\pi_t^k N^{kT}[1+O(1/N)]}{N^{kT}[1+O(1/N)]} = \pi_t^k\left[1+O\left(\frac{1}{N}\right)\right]. \qquad (4)$$

## 5. Conclusion

From (4) we conclude that for any fixed $k$,

$$P_k = \mathbb{P}(I_1 = \cdots = I_k = 1) = \pi_t^k + O\left(\frac{1}{N}\right) \quad (N \to \infty).$$

In particular, for $k = 2$,

$$\mathrm{Cov}(I_i, I_j) = P_2 - \pi_t^2 = O\left(\frac{1}{N}\right).$$

Moreover, by an analogous argument applied to all $2^k$ possible patterns (e.g., using the inclusion–exclusion principle or a straightforward extension), the finite-dimensional distribution of $(I_1, \ldots, I_k)$ converges to the product of independent Bernoulli$(\pi_t)$ distributions. Hence the indicators are asymptotically independent and identically distributed.

**Remark.** The proof relies on the assumptions that $T, t, k$ are fixed and that $p$ is fixed in $(0, 1)$ (so that $1/p$ and $1/(1-p)$ are bounded). If $p$ varies with $N$, we need to assume that $p$ stays bounded away from $0$ and $1$ in order to maintain the same $O(1/N)$ error bound.

Supplementary Table S1: Moment comparison for $X_{=t}$ and $X_{\geq t}$ with equal subset sizes ($N = 100, m = 20, T = 5$).

| t | Method | $X_{=t}$ mean $\mu_t$ | $X_{=t}$ var $\sigma_t^2$ | $X_{=t}$ 3rd c.m. | $X_{=t}$ 4th c.m. | $X_{\geq t}$ mean | $X_{\geq t}$ var | $X_{\geq t}$ 3rd c.m. | $X_{\geq t}$ 4th c.m. |
|---|---|---|---|---|---|---|---|---|---|
| 2 | 公式 | 20.4800000000 | 11.0490232739 | 0.31470596910 | 359.8852999539 | 26.2720000000 | 5.95190051160 | -0.02325724870 | 105.5496648445 |
| 2 | 精确 | 20.4800000000 | 11.0490232739 | 0.31470596910 | 359.8852999539 | 26.2720000000 | 5.95190051160 | -0.02325724870 | 105.5496648446 |
| 2 | 模拟 | 20.4846000000 | 11.1145739857 | 0.65019183540 | 363.0905663248 | 26.2729100000 | 5.94976962960 | 0.18575752290 | 105.2627806518 |
| 3 | 公式 | 5.12000000000 | 3.55746405760 | 1.55053057950 | 37.5680356004 | 5.79200000000 | 3.59053912700 | 1.22510223980 | 38.1763665922 |
| 3 | 精确 | 5.12000000000 | 3.55746405760 | 1.55053057950 | 37.5680356004 | 5.79200000000 | 3.59053912700 | 1.22510223980 | 38.1763665922 |
| 3 | 模拟 | 5.11718000000 | 3.57720461960 | 1.58517304690 | 37.8082772413 | 5.78831000000 | 3.61669351080 | 1.25720697060 | 38.4881022695 |
| 4 | 公式 | 0.6400 | 0.5914 | 0.5038 | 1.4031 | 0.6720 | 0.6170 | 0.5186 | 1.494773 |

| $t$ | Method | $X_{=t}$ mean $\mu_t$ | $X_{=t}$ var $\sigma_t^2$ | $X_{=t}$ 3rd c.m. | $X_{=t}$ 4th c.m. | $X_{\geq t}$ mean | $X_{\geq t}$ var | $X_{\geq t}$ 3rd c.m. | $X_{\geq t}$ 4th c.m. |
|---|---|---|---|---|---|---|---|---|---|
|  |  | 000000000 | 0.622301 | 0.50574460 | 1.4069734 3 | 00000000 | 0.3378877 | 0.7404100 | 5616 |
| 4 | 精确 | 0.6400000000 | 0.59146223001 | 0.503805744600 | 1.40316973430 | 0.67200000000 | 0.61703378877 | 0.51867404100 | 1.4947735616 |
| 4 | 模拟 | 0.63844000000 | 0.59305802972 | 0.51338106380 | 1.44300290460 | 0.67113000000 | 0.61974072050 | 0.52995353695 | 1.5415456753 |
| 5 | 公式 | 0.03200000000 | 0.03180085420 | 0.03140581070 | 0.03365926180 | 0.03200000000 | 0.03180085420 | 0.03140581070 | 0.0336592618 |
| 5 | 精确 | 0.03200000000 | 0.03180085420 | 0.03140581070 | 0.03365926180 | 0.03200000000 | 0.03180085420 | 0.03140581070 | 0.0336592618 |
| 5 | 模拟 | 0.03268900000 | 0.03238168770 | 0.03175942590 | 0.03364837210 | 0.03268900000 | 0.03238168770 | 0.03175942590 | 0.0336483721 |

| $t$ | Method | $X_{=t}$ mean $\mu_t$ | $X_{=t}$ var $\sigma_t^2$ | $X_{=t}$ 3rd c.m. | $X_{=t}$ 4th c.m. | $X_{\geq t}$ mean | $X_{\geq t}$ var | $X_{\geq t}$ 3rd c.m. | $X_{\geq t}$ 4th c.m. |
|---|---|---|---|---|---|---|---|---|---|

Supplementary Table S2: Moment comparison for $X_{=t}$ and $X_{\geq t}$ with equal subset sizes ($N = 5000, m = 1000, T = 5$).

| $t$ | Method | $X_{=t}$ mean $\mu_t$ | $X_{=t}$ var $\sigma_t^2$ | $X_{=t}$ 3rd c.m. | $X_{=t}$ 4th c.m. | $X_{\geq t}$ mean | $X_{\geq t}$ var | $X_{\geq t}$ 3rd c.m. | $X_{\geq t}$ 4th c.m. |
|---|---|---|---|---|---|---|---|---|---|
| 2 | 公式 | 1024.0000000000 | 552.147346 7178 | 16.4307229 519 | 914283.45800 78125 | 1313.6000000000 | 297.406555 5912 | -1.0309391 022 | 265315.67382 81250 |
| 2 | 模拟 | 1024.0981100000 | 553.327477 7027 | 143.162437 7318 | 917842.86128 86615 | 1313.6679000000 | 296.031629 9063 | -3.2304290 422 | 263893.10475 40406 |
| 3 | 公式 | 256.000000000 | 177.367041 7211 | 76.7923042 104 | 94357.027935 0281 | 289.600000000 | 178.475765 8989 | 60.496445 5515 | 95536.469913 4827 |
| 3 | 模拟 | 255.987150000 | 177.247477 3523 | 74.0717540 448 | 94631.814474 6072 | 289.569790000 | 179.044379 7997 | 63.952410 6443 | 96433.903014 9324 |
| 4 | 公式 | 32.00000000 | 29.4928383 999 | 24.9890344 402 | 2626.8191427 379 | 33.60000000 | 30.7547503 885 | 25.695331 9595 | 2854.7953981 399 |
| 4 | 模拟 | 31.98581000 00 | 29.6879055 230 | 25.8470343 572 | 2686.6089809 430 | 33.58264000 00 | 30.9081997 124 | 27.218071 3758 | 2910.0439192 337 |
| 5 | 公式 | 1.60000000 00 | 1.58926437 42 | 1.56798930 46 | 9.1033052716 | 1.60000000 00 | 1.58926437 42 | 1.5679893 046 | 9.1033052716 |

| $t$ | Method | $X_{=t}$ mean $\mu_t$ | $X_{=t}$ var $\sigma_t^2$ | $X_{=t}$ 3rd c.m. | $X_{=t}$ 4th c.m. | $X_{\geq t}$ mean | $X_{\geq t}$ var | $X_{\geq t}$ 3rd c.m. | $X_{\geq t}$ 4th c.m. |
|---|---|---|---|---|---|---|---|---|---|
| 5 | 模拟 | 1.59683000 00 | 1.58159976 71 | 1.53612529 52 | 8.9732 659896 | 1.59683000 00 | 1.58159976 71 | 1.5361252 952 | 8.9732 659896 |

| $t$ | Method | $X_{=t}$ mean $\mu_t$ | $X_{=t}$ var $\sigma_t^2$ | $X_{=t}$ 3rd c.m. | $X_{=t}$ 4th c.m. | $X_{\geq t}$ mean | $X_{\geq t}$ var | $X_{\geq t}$ 3rd c.m. | $X_{\geq t}$ 4th c.m. |
|---|---|---|---|---|---|---|---|---|---|

**Supplementary Table S3:** Moments of $X_{=t}$ and $X_{\geq t}$ for unequal subset sizes ($N = 100$, $\mathbf{m} = [5,20,40,70,30]$, $T = 5$).

| $t$ | $X_{=t}$ mean $\mu_t$ | $X_{=t}$ var $\sigma_t^2$ | $X_{=t}$ 3rd c.m. | $X_{=t}$ 4th c.m. | $X_{\geq t}$ mean | $X_{\geq t}$ var | $X_{\geq t}$ 3rd c.m. | $X_{\geq t}$ 4th c.m. |
|---|---|---|---|---|---|---|---|---|
| 2 | 37.270000000 | 20.9816510000 | 2.59253200000 | 1309.6921300000 | 54.694000000 | 7.159947000 | -0.15934800000 | 152.9662450000 |
| 3 | 15.050000000 | 7.9082560000 | 0.86444400000 | 184.9611920000 | 17.424000000 | 6.314953000 | 0.43902200000 | 118.7886720000 |
| 4 | 2.2900000000 | 1.8553500000 | 1.19113400000 | 10.6441570000 | 2.3740000000 | 1.888293000 | 1.16802500000 | 10.9604050000 |
| 5 | 0.0840000000 | 0.0821300000 | 0.07848900000 | 0.0917290000 | 0.0840000000 | 0.082130000 | 0.07848900000 | 0.0917290000 |